\documentclass{article}
\usepackage{graphics}
\usepackage{amsmath}
\usepackage{amssymb}
\usepackage{amsfonts}
\newcommand{\bea} {\begin{eqnarray*}}
\newcommand{\beq} {\begin{equation}}
\newcommand{\bey} {\begin{eqnarray}}
\newcommand{\eea} {\end{eqnarray*}}
\newcommand{\eeq} {\end{equation}}
\newcommand{\eey} {\end{eqnarray}}

\newcommand{\raz}{{\rm 1\kern-0.25em I}}

\begin{document}
\title {Improved Lower Bounds for the Critical
Probability of Oriented-Bond Percolation in Two Dimensions }
\author{V. Belitsky\thanks{Supported by CNPq (grant N.301637/91-1) / e-mail: belitsky@ime.usp.br}\\
Universidade de S\~ao Paulo - Instituto de Matem\'atica e Estat\'\i stica\\
\and  T. L. Ritchie\thanks{Supported by a grant from CNPq / e-mail: thomasl@ime.usp.br}\\
Universidade de S\~ao Paulo - Instituto de Matem\'atica e Estat\'\i stica\\
\small Rua do Mat\~ao 1010, CEP 05508-900, S\~ao Paulo SP, Brazil. }
\normalsize
\date{} 
\maketitle \begin{abstract}
  We present a coupled decreasing sequence of random walks on $ \mathbb Z $
  ($\overline{X}^{(i)}_{.} ,i \in \mathbb N  $)  that dominates
the edge process of oriented-bond percolation in two dimensions.
  Using the concept of \textquotedblleft random walk in a strip \textquotedblright ,we construct an algorithm
  that generates an increasing sequence of lower bounds that converges to the
  critical probability of oriented-bond percolation.
  Numerical calculations of the first ten lower bounds thereby generated lead
  to an improved,ie higher, rigorous lower bound to this critical probability, viz. $p_{
 c} \geq 0.63328 $. Finally a computer simulation technique is presented; the
  use thereof establishes $0.64450$ as a non-rigorous five-digit-precision (lower) estimate for $p_{c}$.
\newline
\newline
\underline{$\ $ \hspace{140mm}$\ $}\newline
\underline{Markov chain in a strip, coupling, simulation .\hspace{77.5mm}}
\end{abstract}

\section{Introduction}

Oriented percolation in two dimensions or discrete time contact process on
  $ \mathbb Z $ is a one-parameter family of discrete time stochastic
  processes defined on $ \{0,1\}^{\mathbb Z} $; the parameter $p$, the
  infection rate, taking values on $[0,1]$. It is a well established fact that
  this family exhibits a fase transition, as the value of $p$ increases from
  $0$ to $1$: if $p<p_{c}\ \ $\footnote{``$p_{c}$'' is known as the critical
  probability of the family.  }, the process dies out almost surely; whereas if
  $p>p_{c}$, the process survives with positive probability (see
  ~\cite{Durrett_1} for instance ).

As usual, in critical phenomena theory, an analytical expression for $p_{c}$
is unknown and its value has been estimated both in mathematical and 
physical literatures ~\cite[Sec.6]{Durrett_1}.

Up to the present moment, the  best (rigorous) lower and upper bounds for
$p_{c}$ are respectively $0.6298$ and $2/3$ ~\cite{Grimmett_1}.

In this paper, we present an algorithm that generates an increasing sequence
of lower bounds that converges to the critical probability of oriented bond
percolation in two dimensions and, calculating the first ten lower bounds
thereby generated, we were able to improve the best lower bound known up to
now from $0.6298$ to $0.63328$.

More specifically:
\begin{enumerate}
\item a numerical sequence of lower bounds for $p_{c}, \{ p_{c}^{(i)}\}_{i\in\mathbb N}$,  was constructed in the following steps:
      \begin{description}
      \item[(a)] Associated to a decreasing\footnote{in a sense to be
          specified later by means of coupling} sequence(in $i$)  of random walks on $\mathbb
        Z$, $ \{\overline{X}^{(i)}_{n}\}_{n\in\mathbb N}$, that
        dominate the edge process of oriented percolation, we constructed a
        sequence (in $i$) of finite, discrete time
        Markov chains, $\{Y^{(i)}_{n}\}_{n\in\mathbb N}$, 
        whose transition probabilities are specific\footnote{polinomial}
        functions of $p$, the infection probability; and
      \item[(b)]calculated the mathematical expectation of the random variable \linebreak
        $\mathbb
        E\left(\overline{X}^{(i)}_{1}-\overline{X}^{(i)}_{0}|Y^{(i)}_{0}\right)$, the
        mean jump on the $Y$-configuration,
        with respect to $\pi^{(i)}$, the stationary measure for $Y^{(i)}_{\bullet}$.
      \item[(c)] $p^{(i)}_{c}$ was then defined to be the (only) value of $p$
        that nullifies the above expectation. In other words, when
        $p=p^{(i)}_{c}$, the random walk
        $\{\overline{X}^{(i)}_{n}\}_{n\in\mathbb N}$ has zero mean drift.
      \end{description}
\item the numerical sequence $\{p^{(i)}_{c}\}_{i\in\mathbb N}$ was shown to
  converge in a non-decreasing fashion to $p_{c}$,\linebreak ie $p^{(i)}_{c}\nearrow
  p_{c}$. Moreover,
\item the first ten lower bounds $(p_{c}^{(0)},p_{c}^{(1)},...,p_{c}^{(9)})$ were
  numerically calculated, thereby improving the best rigorous lower bound 
  known up to the moment from $0.6298$ to $0.63328$.
\end{enumerate}
In the last section we present a simulation technique which in our opinion has
important advantages in comparison with the usual so called \emph{Monte Carlo}
simulation techniques, in the sense that it exhibites a clear cut off between
the subcritical and supercritical phases; thereby enabling a precise
estimation of the critical probability of \emph{ Oriented Bond Percolation}
without the aid of scaling techniques. By means thereof a lower bound for
$p_{c}$ was obtained within a precision of 5 digits, \emph{viz} $p_{c}^{1000}=0.64451$, so that we can (non-rigorously) state that $p_{c}=(0.64451\pm 0.00001)$.

\section{Definitions and Constructions}
\subsection{The Enviroment}
\label{Omeio}
Let $\mathcal{G}\equiv(\mathcal{V},\mathcal{E})$ be the oriented graph, having
\(\mathcal{V}=\{(n,m):n\in\mathbb N \mbox{ and $(n+m)$ is even}\}\)  as its set
of vertices/sites, and  \(\mathcal{E}=\{ e^{l}_{nm},e^{r}_{nm} : (n,m) \in
\mathcal{V} \}\)  as its set of bonds.

Bond $e^{l}_{n,m}$ points from site $(n,m)$ to site $(n+1,m-1)$, whereas bond
$e^{r}_{n,m}$ points from site $(n,m)$ to site $(n+1,m+1)$. Sometimes the
\textit{natural} association: $l\leftrightarrow -1\ /\ r\leftrightarrow+1 $ will be
assumed through out this text.

It is useful to think of $n$ as a (discrete) time coordinate and of $m$
as a  (discrete) space coordinate of the graph $\mathcal{G}$.

 $\mathcal{V}_{n}$ denotes the n-th slice of $\mathcal{V}$, ie  $\mathcal{V}_{n}=
 \{(n,m)\in\mathcal{V}:n \ is fixed\: \}$, and
$\mathbb Z_{n}$ the set of integers $m$ such that  $(m+n)$ is even.

It will be often useful in the forthcoming definitions to identify
$\mathcal{V}_{n}$ as $\mathbb Z_{n}$, and to think of
$\{0,1\}^{\mathcal{V}_{n}}$ as being $\{0,1\}^{\mathbb Z_{n}}$. Therefore
arises the loose, yet natural, notation:\vspace{1 mm}\linebreak
\hspace{-4 mm} $\eta_{n}(m)\equiv\eta_{n}(n,m),\ \eta_{n}\in\{0,1\}^{\mathcal{V}_{n}}$.
\vspace{3.3cm}
\pagebreak 
\begin{figure}[h]
\hspace{5cm}
\scalebox{.3}{\includegraphics[3.9cm,0.9cm][30.6cm,20.1cm]{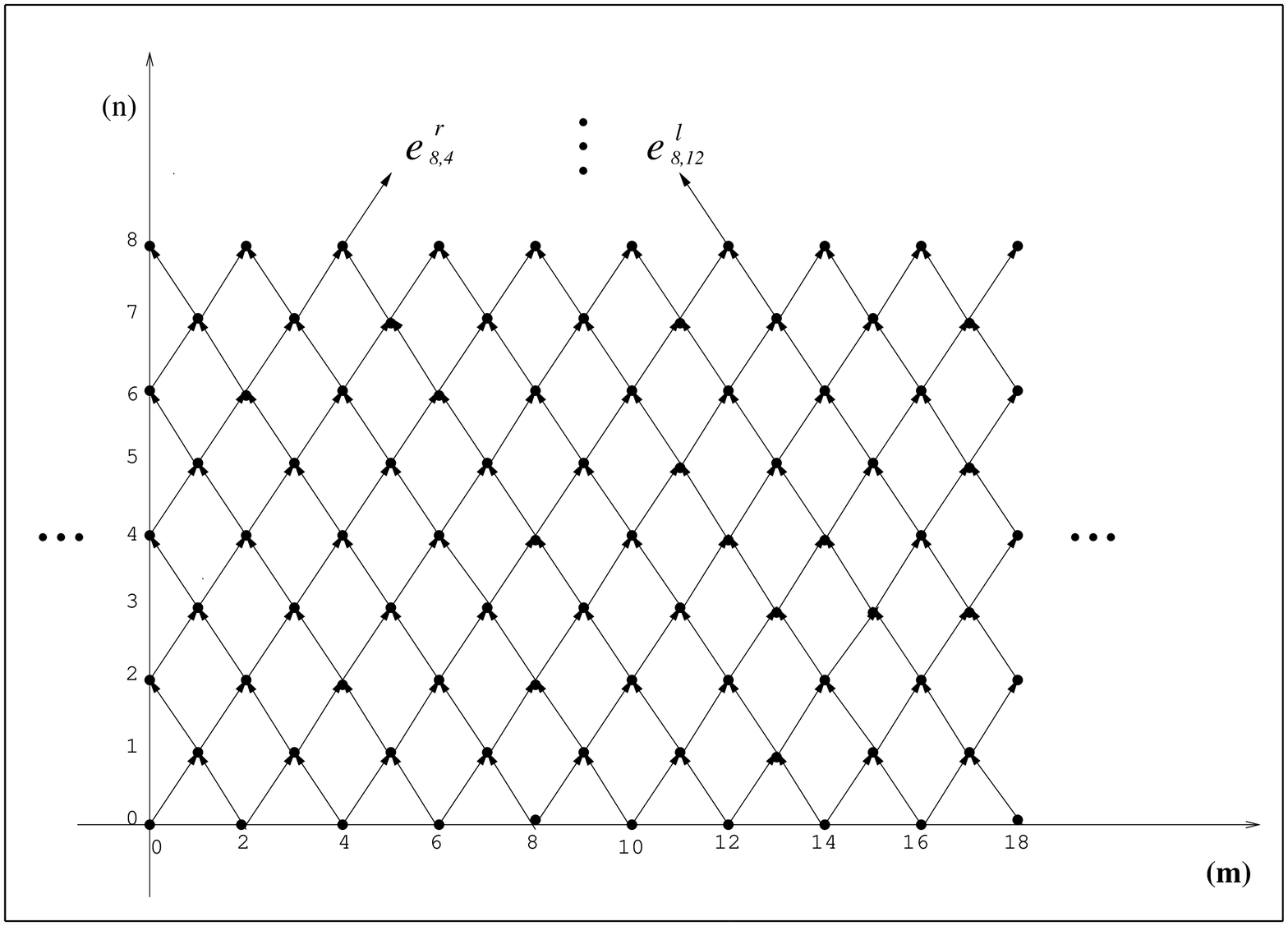}}
\small 
\caption{ Graph $\mathcal{G}$, whereon the processes  will be defined
  (Sec.~\ref{processos} below).}
\normalsize
\end{figure} 

\subsection{The Probability Structure}
 Let \( \{ \xi^{j}_{nm}:(n,m)\in \mathcal{V},j \in \{l,r\} \}\) be a family of
 independent and uniformly distributed \linebreak  ( onto $[0,1]$ ) random variables,
 defined on the same abstract probability space $\Omega$. Starting from this
 three-indexed countable family, we construct the four-indexed uncountable
 family of \emph{iid Bernoulli} random variables on the same probability space
 $\Omega$ :
\begin{equation}
\{ _{p}\xi^{j}_{n,m}:(n,m)\in\mathcal{V},j\in \{l,r\},p \in [0,1] \}
\end{equation}
such that
\begin{equation}
 _{p}\xi^{j}_{n,m}=\mbox{\large\textbf{1}}_{ \{ \xi^{j}_{nm}\leq p\} }
\label{indicador}
\end{equation}  
\newline
It follows straightfowardly from (~\ref{indicador})  that

\begin{equation}
\left\{\begin{array}{lr}         
\mathbb P (\mbox{ }_{p}\xi^{j}_{nm} = 1)\mbox{ }=\mbox{ }p & (a)\\ 
\mathbb P (\mbox{ }_{p}\xi^{j}_{nm} = 0)\mbox{ }=\mbox{ }1-p\mbox{ }:=\mbox{
  }q & (b)\end{array}\right. 
\end{equation}
and the \textbf{Fundamental Coupling:}
\newline
\begin{equation}
p_{1} \leq p_{2} \Rightarrow \mbox{ }_{p_{1}}\xi^{j}_{nm}(\omega) \leq  \mbox{
  }_{p_{2}}\xi^{j}_{nm}(\omega),\mbox{ }\forall\omega\in\Omega,\mbox{
  }(n,m)\in\mathcal{V},\mbox{ }j\in\{l,r\}
\end{equation}
\newline
  We observe that $\mbox{ }_{p}\xi^{l}_{nm}=1 $ is interpreted as a channel
  open to infection propagation from site $(n,m)$ to site 
  $(n+1,m-1)$; $\mbox{ }_{p}\xi^{r}_{nm}=0$, as a channel obstructed to
  infection propagation from site $(n,m)$ to site $(n+1,m+1)$; and so forth ...  
\subsection{The Processes}
\label{processos}
\subsubsection{The Strengthened Discrete Time Contact Processes (SDTCP)}

$\eta_{n}\in \{0,1\}^{\mathcal{V}_{n}}$, is interpreted as a \emph{infection}
state on slice \(\mathcal{V}_{n}\):
\begin{equation}
\left\{ \begin{array}{l}
        \eta_{n}(m)=1 \mbox{:  site \emph{m} is infected at time \emph{n}}
      \\\eta_{n}(m)=0 \mbox{;  site \emph{m} is healthy at time \emph{n}}
        \end{array}
\right.  
\end{equation}  

Particulary $\eta_{0}\in\{0,1\}^{\mathcal{V}_{0}}$ will denote an initial state
of infection over $\mathcal{V}_{0}$; the set of even integers, according to
Section ~\ref{Omeio} above.

At this point, the following notation (to be used throughout this paper)
should be kept into account: $\mathbb N =\{0,1,2,\cdots \}$ and
$\overline{\mathbb N}= \mathbb N \cup \{\infty \}$. 

Now for each $i\in\ \overline{\mathbb N}$, we define (by induction on $n$) the
sequence of $\{0,1\}^{\mathcal{V}_{n}}$-valued random variables
$\{\:_{p}X^{(i)}_{n}\}_{n\in\mathbb N\ }$ as
\newtheorem{df}{Definition}[section]
\begin{df}
\label{dfrec}
.

\begin{description}
\item[(a)]$ \mbox{  }_{p}X^{(i)}_{0}:=\eta_{0}$, so that   $   \mbox{
    }_{p}X^{(i)}_{0}(m)=\eta_{0}(m)\mbox{ },\forall m\in \mathbb Z_{0}$
\item[(b)]$\mbox{  }_{p}X^{(i)}_{n+1}(m)=\sup \left\{ \mbox{
    }_{p}X^{(i)}_{n}(m-1).\mbox{ }_{p}\xi^{r}_{n,m-1};\mbox{
    }_{p}X^{(i)}_{n}(m+1).\mbox{ }_{p}\xi^{l}_{n,m+1};\mbox{ }\mbox{\Large
\textbf{1}}_{\{\mbox{ }_{p}\overline{X}^{(i)}_{n+1}-m>2i\}} \right\}$, \newline where
  $\mbox{ }_{p}\overline{X}^{(i)}_{n+1}=\sup_{m\in\mathbb
    Z_{n};j\in\{-1,+1\}} \left\{ (m+j):\mbox{
      }_{p}X^{(i)}_{n}(m).\xi^{j}_{nm}=1 \right\}$
\end{description}
\end{df}
The role of the indicator function in Definition ~\ref{dfrec}(b) above is to
 infect \emph{by force} all the sites lying farther than $2i$ on the left side
 of $\:_{p}\overline{X}^{(i)}_{n+1}$, the utmost right infected site at time
 $n+1$. It is natural,therefore, to call the stochatic process
 $ \{\ _{p}X^{(i)}_{n}\}_{n\in \mathbb N}$, defined above, as the
\textbf{Strengthened Discrete Time Contact Process} (of i-th order and
infection parameter $p$ ).

It follows directly from Definition~\ref{dfrec} that
\begin{description}
\item[(a)]
    when $i=\infty$, the indicator function does not act anymore, and we recover the ordinary \textbf{Discrete Time
    Contact Process}, also called \textbf{Oriented Percolation in Two
    Dimensions}, described in  ~\cite[Sec. 2]{Durrett_1}. So forth in this paper we shall
    refer to it as $\{\:_{p}X^{(\infty)}_{n}\}_{n\in \mathbb N}$ ;
\item[(b)]
 the sequence of SDTCPs $\{\:_{p}X^{(i)}_{n}\}_{n \in \mathbb N}$ is
  decreasing in $i$, in the following sense:\small
\begin{equation}
\label{acopl1}
i\leq j;i,j\in \overline{\mathbb  N} \Rightarrow  \mbox{ }_{p}X^{(i)}_{n}(m)[\omega]\geq \mbox{
  }_{p}X^{(j)}_{n}(m)[\omega],\forall \omega \in \Omega,(n,m)
  \in\mathcal{V},p\in[0,1] 
\end{equation}\normalsize
In particular, the \textit{ Oriented Percolation Process} is the weakest of
them all.

 Inequality ~\ref{acopl1} above is called \textbf{Coupling of First Kind}
\item[(c)]
the family of SDTCPs $\{\:_{p}X^{(i)}_{n}\}_{n\in {\mathbb N}}$ is increasing in
$p$ in the following sense
\begin{equation}
\label{acopl2}
p_{1} \leq p_{2} \Rightarrow \:_{p_{1}}X^{(i)}_{n}(m)[\omega] \leq
\:_{p_{2}}X^{(i)}_{n}(m)[\omega],\forall \omega \in \Omega,(n,m) \in
\mathcal{V},i \in \overline{\mathbb N}
\end{equation}

 Inequality ~\ref{acopl2} above is called \textbf{Coupling of Second Kind}
\end{description}
In this paper, unless otherwise stated, $\eta_{0}$ will be chosen to be
$\mbox{\large\textbf{1}}_{\bullet\leq 0}\in\{0,1\}^{\mathcal{V}_{0}}$, where
\begin{equation}
\label{meiocheio}
\mbox{\large\textbf{1}}_{\bullet\leq 0}(m)= \left\{ \begin{array}{l}
1,m\leq0\\0,m>0 \end{array} \right.
\end{equation}
We shall write $\{\:^{\eta}_{p}X^{(i)}_{n}\}_{n \in \mathbb N}$ to emphasize
that $\eta_{0}=\eta$, when $\eta$ is some specific element of
$\{0,1\}^{\mathcal{V}_{0}}$ ( probably different from
$\mathbf{1}_{\bullet\leq 0}$). Therefore the symbols
$\{\:_{p}X^{(i)}_{n}\}_{n \in \mathbb N}$,
$\{\:^{\mathbf{1}_{\bullet\leq 0}}_{\ \ \ \ p}X^{(i)}_{n}\}_{n \in \mathbb
  N}$, $\:^{\mathbf{1}_{\bullet\leq 0}}_{\ \ \ \ p}X^{(i)}_{\bullet}$,
$\:_{p}X^{(i)}_{\bullet}$ have all the same meaning and, for simplicity's
sake, the last will be chosen, when no confusion may arise.

As usual in Particle System's notation we write $\eta_{n}\geq\theta_{n}$, whenever
$\eta_{n}(m)\geq\theta_{n}(m)$ for every $m\in\mathbb Z_{n}$.

Definition ~\ref{dfrec} also implies that, for $\eta,\theta\in\{0,1\}^{\mathcal{V}_{0}}$ 
\begin{equation}
\label{acopl3}
\eta\geq\theta\Rightarrow\:^{\eta}_{p}X^{(i)}_{n}(m)[\omega]\geq
\:^{\theta}_{p}X^{(i)}_{n}(m)[\omega],\forall
p\in[0,1],(n,m)\in\mathcal{V},\omega\in\Omega,i\in\overline{\mathbb N}
\end{equation}

Inequality ~\ref{acopl3} is called \textbf{Coupling of Third Kind}.
\begin{figure}

\hspace{3cm}
\scalebox{.4}{\includegraphics[25.6cm,16.5cm]{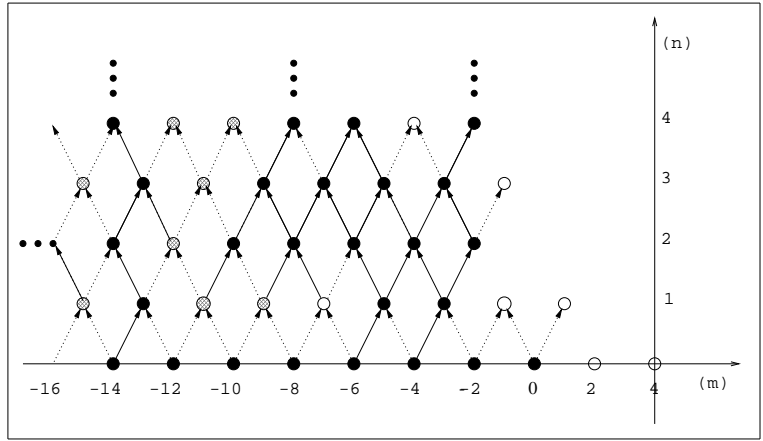}}
\caption{\label{figura2} \small A realization $ \omega $ of the $\:_{p}X^{(2)}_{\bullet}$  process . Full-black arrows are open, whereas dotted ones are
  closed to infection propagation. Black sites/cells are infected, white ones
  are healthy, while the \emph{gridded} ones were reinfected.  \normalsize}
\end{figure}

\subsubsection{The Right Edge Processes (REP)}
Given a particular SDTCP $\:_{p}X^{(i)}_{\bullet}$ and assuming that
$\:_{p}\overline{X}^{(i)}_{0} \stackrel{\mathrm{def}}{=}0$, the second line of
Definition ~\ref{dfrec}(b) above defines a (non-markovian) random process on
$\mathbb Z$ denoted by $\{\:_{p}\overline{X}^{(i)}_{n}\}_{n\in\mathbb N}$ ( or
$\:_{p}\overline{X}^{(i)}_{\bullet}$ in abbreviated fashion). In Figure
~\ref{figura2}  above $ \:_{p}\overline{X}^{(2)}_{0}=0 ,
\:_{p}\overline{X}^{(2)}_{1}=-3 ,  \:_{p}\overline{X}^{(2)}_{2}=-2,
\:_{p}\overline{X}^{(2)}_{3}=-3,  \:_{p}\overline{X}^{(2)}_{4}=-2,$ \ldots
\vspace{1 mm}  

Again, in  case of $i=\infty$, the \emph{Edge Process of Oriented Percolation}
cited in ~\cite{Durrett_1} is recovered. Throughout this text, this term will
be employed in a generilized way ($i \in \overline{\mathbb N}$).

It is useful to think of $\:_{p}\overline{X}^{(i)}_{\bullet}$ as random walks
on $\mathbb Z$.
\subsubsection{The Induced Markov Chains}
\label{CMIs} 
\begin{df}
\label{IMC}
The Markov chain $\{\:_{p}Y^{(i)}_{n}\}_{n\in\mathbb N}$, $i\in\mathbb{N}$, with state space
$\mathcal{S}^{(i)}=\{0,1\}^{\{2,4,\ldots,2i\}}$ defined by
\[\:_{p}Y^{(i)}_{n}(2j)=\:_{p}X^{(i)}_{n} \left(
  \:_{p}\overline{X}^{(i)}_{n}-2j \right),1\leq
j\leq i\]
is called \textbf{Induced Markov Chain}(IMC) for the SDTCP $\:_{p}X^{(i)}_{\bullet}$. 
\end{df}
In what follows, a generic element $\sigma\in\mathcal{S}^{(i)}$, will be
labeled by $n\in\{0,1,\ldots,2^{i}-1\}$ according to
\begin{equation}
\sigma=\sigma_{n}\Leftrightarrow n=\sum_{j=1}^{i}2^{i-j}\sigma(2j)
\end{equation}
Accordingly, in Figure ~\ref{figura2} above
$\:_{p}Y^{(2)}_{0}=\sigma_{3},\:_{p}Y^{(2)}_{1}=\sigma_{2},\:_{p}Y^{(2)}_{2}=\sigma_{3},\:_{p}Y^{(2)}_{3}=\sigma_{3},\:_{p}Y^{(2)}_{4}=\sigma_{1},$
\ldots \vspace{1 mm}\newline
The transition probabilities for the IMC $\:_{p}Y^{(i)}_{\bullet}$
\begin{equation}
\label{trans.prob} 
q^{(i)}_{lm}=\mathbb
P \left[ \:_{p}Y^{(i)}_{n+1}=\sigma_{m}\mid\:_{p}Y^{(i)}_{n}=\sigma_{l}\right]
,0\leq
l,m\leq2^{i}-1
\end{equation} 
are all strictly positive polinomial functions of $p$ (provided
$0<p<1$). Therefore $\pi^{(i)}$, its stationary measure on
$\mathcal{S}^{(i)}$, 
is well defined.

The notation
\[\pi^{(i)}_{l}\stackrel{\mathrm{def}}{=}\pi^{(i)}(\sigma_{l}),l=0,1,\ldots,2^{i}-1\]
is self explanatory.

At this point, in accoradance with ~\cite[Sec.3.1]{Misha}, we state 
\begin{df}
\label{Pilmk}
\[P^{(i)}_{(l,m,k)}(p)\stackrel{\mathrm{def}}{=}\mathbb P \left[
  \:_{p}\overline{X}^{(i)}_{n+1}=\:_{p}\overline{X}^{(i)}_{n}+
  (1-2k),\:_{p}Y^{(i)}_{n+1}=\sigma_{m}\mid \:_{p}Y^{(i)}_{n}=\sigma_{l}  \right],\ k\in
  \mathbb N,0\leq l,m \leq 2^{i}-1\]
the transition probabilities from state $\sigma_{l}$ to state $\sigma_{m}$
  with a jump of magnitude $(1-2k)$ 
\end{df}
\begin{df}
\label{Mil}
\[M^{(i)}_{l}(p)\stackrel{\mathrm{def}}{=}\sum^{\infty}_{k=0}
  \sum^{2^{i}-1}_{m=0}(1-2k).P^{(i)}_{(l,m,k)}(p),\ \ 0\leq l \leq2^{i}-1 \]
the mean jump/drift of the REP on configuration $\sigma_{l}$
\end{df}
and
\begin{df}
\label{Mdrift}
\[M^{(i)}(p)\stackrel{\mathrm{def}}{=}\sum^{2^{i}-1}_{l=0}M^{(i)}_{l}(p).\pi^{(i)}_{l}\]
  the mean jump/drift of the REP.
\end{df}
$M^{(i)}(p)$ is also called \emph{the (right)edge speed} in accordance with
Section ~\ref{probcrit} bellow.
As it will soon become clear, Definition ~\ref{Mdrift} is of fundamental
importance in this paper.

Under this framework, the SDTCPs $\:_{p}X^{(i)}_{\bullet}$, described
above, can be regarded, in case of $i \in \mathbb N $, as \emph{Markov
  Chains/Random Walks in a $2^{i}$-rowed Strip}:
\begin{equation}
\:_{p}X^{(i)}_{n}=\:_{p}X^{(i)}_{n}\left(\:_{p}\overline{X}^{(i)}_{n},
  \:_{p}Y^{(i)}_{n}\right)
\end{equation}
a slightly different idea of \emph{Markov Chains in a Half Strip}, developed in
~\cite[Sec.3.1]{Misha}.
\subsection{Critical Probabilities}
\label{probcrit}
In Section ~\ref{Pteorica}  below it will be shown that, for $i\in \mathbb N $ and $p \in
[0,1]$:
\begin{description}
\item[(i)] $M^{(i)}(p)$ is a scrictly incresing function of $p$;
\item[(ii)]$M^{(i)}(p)$ has only one real root ( in $(0,1]$);
\item[(iii)] \Large\(
  \frac{\:_{p}\overline{X}^{(i)}_{n}}{n}\)\normalsize \( \stackrel{(n)}{\longrightarrow}M^{(i)}(p)\ \ \  \
  \) \emph{a s}  
\end{description}

In case of $i= \infty$, $\alpha(p)$, the (right) edge speed of oriented percolation, plays precisely
the same role of $M^{(i)}(p)$ in the finite case depicted
above\footnote{$\alpha(p)=- \infty$,when $p<p_{c}$; so that the strict
  increasing behaviour does not apply to $ \alpha(p)$ precisely. }.(For
details, see ~\cite{Durrett_1}, for instance).So the notation 
\[ M^{(\infty)}(p)\stackrel{\mathrm{def}}{=}\alpha (p)\]
suggests itself and we state the following 
\begin{df}
\label{dfpc}
The critical (infection) probability $ \mathit{p^{(i)}_{c}}$ for the family of
stochastic processes $ \:_{p}X^{(i)}_{\bullet},i \in \overline{\mathbb N}$, is
the only real root of $M^{(i)}(p)$ (in $(0,1]$ ), the edge speed. Hence
\[M^{(i)} \left( p^{(i)}_{c} \right) = 0,\: \forall i \in \overline{ \mathbb
  N}\]
 
\end{df}
The heuristic meaning of Definition ~\ref{dfpc} above is :
\begin{description}
\item[(i)]for $p<p_{c}^{(i)}$, $ \lim_{n \rightarrow \infty}
  \:_{p}\overline{X}^{(i)}_{n} = - \infty \ \ \  \mathit{a\:s} $, ie the infection
  dies out with probability one;
\item[(ii)] for $p>p_{c}^{(i)}$, $ \lim_{n \rightarrow \infty}
  \:_{p}\overline{X}^{(i)}_{n} = + \infty \ \ \ \mathit{a\:s} $, ie the infection
  spreads out over all $\mathbb Z $.
\end{description}
In the sequel, we prove an important relation concerning the \vspace{1mm}
$p_{c}^{(i)}$s, $i \in \overline{\mathbb N}$,  just defined, viz. \mbox{$p_{c}^{(i)}
\nearrow_{i} \  p^{(\infty)}_{c} \ \ $}\footnote{$p^{(\infty)}_{c}
  \stackrel{\mathrm{def}}{=}p_{c}$, putting the notations of Sections 1 and 2
  into agreement.}.
This non-decreasing convergence to the critical probability of \emph{Oriented
  Percolation} (to be called \emph{The Convergence Theorem}), besides the
  possibility of calculating the $p_{c}^{(i)}$s ($i \in \mathbb N$) by
  algebraical means \mbox{( Section ~\ref{CalcNum} below)}, is the cornerstone of this work.

\section{The Convergence Theorem and Preliminary Results}
\label{Pteorica}
\newtheorem{lem}{Lemma}[section]  
\newtheorem{corolario}[lem]{Corollary}
\newtheorem{teorema}[lem]{Theorem}
\begin{lem}
\label{lema1}
 \Large
 \( \frac{\:_{p}\overline{X}^{(i)}_{n}}{n}\)\normalsize\(
 \stackrel{(n)}{\longrightarrow}M^{(i)}(p)\ \ \   \) \emph{a s}, $\forall i
 \in \overline{\mathbb N},\:p \in (0,1]$  
\end{lem}
\small
\emph{ Proof:} \begin{description}\item[ First case ($i \in \mathbb
  N$): ] Let $n^{(i)}_{l}$ be the (random) number of visits that the IMC
$\:_{p}Y^{(i)}_{\bullet}$ makes to state $\sigma_{l}$ up to time $n$
 ( so that
\( \sum^{2^{i}-1}_{l=0}n^{(i)}_{l}=n \)) and \( J^{(i)}_{l,k} \) the $k^{th}$
jump of the REP $ \: _{p} \overline{X}^{(i)}_{\bullet} $ on state/row
$\sigma_{l}$.\newline
The (strong) Markov property of  SDTCP $\: _{p} X^{(i)}_{\bullet}$ makes  $
J^{(i)}_{l,k},:\ k \in \{1,2,3,\ldots \}$ \emph{iid} rvs with $ \mathbb E
\left[ J^{(i)}_{l,k} \right]=M^{(i)}_{l}$.\newline
Now
\begin{eqnarray*}
  \frac{\:_{p}\overline{X}^{(i)}_{n}}{n} = \sum^{2^{i}-1}_{l=0} 
\sum^{n^{(i)}_{l}}_{k=1}\frac{J^{(i)}_{l,k}}{n} =
 \sum^{2^{i}-1}_{l=0}
 \sum^{n^{(i)}_{l}}_{k=1}\frac{J^{(i)}_{l,k}}{n^{(i)}_{l}}. 
\frac{n^{(i)}_{l}}{n} = \sum^{2^{i}-1}_{l=0}\frac{n^{(i)}_{l}}{n}
 \sum^{n^{(i)}_{l}}_{k=1}\frac{J^{(i)}_{l,k}}{n^{(i)}_{l}}\ , \\ \frac{n^{(i)}_{l}}{n} \stackrel{(n)}{\longrightarrow} \pi^{(i)}_{l} \
\mathit{a\:s\ } \mbox{ and}\  \sum^{n^{(i)}_{l}}_{k=1} \frac{J^{(i)}_{l,k}}{n^{(i)}_{l}}
\stackrel{n^{(i)}_{l}}{\longrightarrow} M^{(i)}_{l}(p) \ \mathit{a\:s\ }
\end{eqnarray*}
Taking the limit $n \rightarrow \infty$ and bearing definition
~\ref{Mdrift} in mind, we get the desired result.
\item[Second case ($i=\infty$): ] See ~\cite[pag.1004]{Durrett_1}
\end{description}
\begin{flushright} $ \Box $ \end{flushright}    
\normalsize
\begin{lem}
\label{lema2}
  For $i \in \mathbb N $, each function $ \begin{array}{rcl}
                \ & \ & \ \\
                M^{(i)}:(0,1] & \rightarrow & (-\infty,1]\\
                    \ \ \ p & \mapsto & M^{(i)}(p)
      
                 \end{array}$ is strictly increasing (in $p$). \vspace{2mm}Moreover,
                $M^{(i)}$ is a surjection from $(0,1]$ onto $(- \infty,1]$.
\end{lem}
\small
\emph{Idea of Proof:}\par The non-decreasing behaviour of $M^{(i)}(p)$ follows
from the \emph{Second Kind of Coupling} (inequality ~\ref{acopl2} above) and
Lemma ~\ref{lema1} just established. Definitions ~\ref{IMC},~\ref{Pilmk},~\ref{Mil} and
~\ref{Mdrift} make $M^{(i)}(p)$ a rational function of $p$ such that
$M^{(i)}(0_{-})=- \infty $
and $M^{(i)}(1)=(1)$; so $M^{(i)}(p)$ can not be constant on any interval
$[p_{1},p_{2}] \subset (0,1] $ and the  strict behaviour follows.
\begin{flushright} $ \Box $ \vspace{-5 mm} \end{flushright}
\normalsize
\textbf{Commentary on Lemma ~\ref{lema2} :} \newline
At this point,it is worth observing that the function $ \alpha(p)
\stackrel{\mathrm{def}}{=} M^{(\infty)}(p) $ is non-decreasing; strictly
positive, when $p>p^{c}$; null, when $p=p_{c}$ and infinetly negative, when
$p<p_{c}$. Again, the reference is ~\cite{Durrett_1}.
  
Lemma~\ref{lema2} above yields
\begin{corolario}
\label{coro3}
 For each $i \in \overline{\mathbb N}$, $M^{(i)}$ has only one real root in (0,1], denoted by
 $p_{c}^{(i)}$,  in
accordace with \mbox{Section ~\ref{probcrit}}.
\end{corolario}
\begin{lem}
\label{lema4}
The sequence of functions $\left\{ M^{(i)}(p) \right\}_{i \in \overline{\mathbb N}} $, is
non increasing, \newline ie \mbox{$ i \leq j;\: i,j \in \overline{\mathbb N}
  \Rightarrow  M^{(i)}(p)\geq M^{(j)}(p),\: \forall p \in (0,1] $.} In
particular $\alpha(p) \leq M^{(i)}(p);\forall i \in \mathbb N,\:\forall p \in (0,1]$.
\end{lem}
\small
\emph{Idea of Proof:}\par The non-increasing behaviour (in $i$) of $\left\{
  M^{(i)}(p) \right\}$ follows from the \emph{First Kind of Coupling} (
inequality~\ref{acopl1} above) and again from Lemma~\ref{lema1}.
\begin{flushright} $\Box$ \end{flushright}
\normalsize
Lemma ~\ref{lema4} yields
\begin{corolario}
\label{coro5}
The numerical sequence $ \left\{ p^{(i)}_{c} \right\}_{i \in \mathbb N}$ is
non-decreasing. Hence $ \: \lim_{i \rightarrow \infty}p^{(i)}_{c}$ is well
defined and $ p^{(i)}_{c} \nearrow_{i} \: \lim_{i \rightarrow
  \infty}p^{(i)}_{c}$. Moreover $p_{c} \stackrel{\mathrel{def}}{=}  p^{(\infty)}_{c} \geq p^{(i)}_{c},\: \forall
 i \in \mathbb N $
 \end{corolario}
and
\begin{corolario}
\label{coro6}
$p^{(\infty)}_{c} \geq \lim_{i \rightarrow \infty}p^{(i)}_{c}.$
\end{corolario}
Now we turn our attention to the reverse (and more difficult) inequality,
viz.$\:p_{c} \leq \lim_{i \rightarrow \infty}p^{(i)}_{c}$. For that
 consider the following probability spaces:
\begin{itemize}
\item $(\Omega,\mathbb P)$ : the abstract probability space, where the $
  \mathit{iid}$ \textit{rv}s $\xi^{j}_{nm}$ were defined;
\item $(\mathcal{S}^{(i)},\pi^{(i)})$ : the finite probability space of the
  IMC $Y^{(i)}_{\bullet}$, endowed with its stationary measure $\pi^{(i)}$;
\item $(\Omega \times \mathcal{S}^{(i)},\mathbb P \times \pi^{(i)})$ : the
  product space.
\end{itemize}
\pagebreak
In the product space, we define the stochastic processes :
\begin{itemize}
\item $\left\{ \:_{p}\overline{\mathcal{X}}^{(i)}_{n} \right\}_{n \in \mathbb
    N}$, by $ \:_{p}\overline{\mathcal{X}}^{(i)}_{n} \left(
    (\omega,\sigma) \right) \stackrel{\mathrm{def}}{=}
  \:_{p}\overline{X}^{(i)}_{n}(\omega),\ \forall(\omega.\sigma) \in \Omega
  \times \mathcal{S}^{(i)}$
\item  $\left\{ \:^{\pi^{(i)}}_{p}\overline{\mathcal{X}}^{(i)}_{n} \right\}_{n \in \mathbb
    N}$, by  $ \:^{\pi^{(i)}}_{p}\overline{\mathcal{X}}^{(i)}_{n} \left(
    (\omega,\sigma) \right) \stackrel{\mathrm{def}}{=}
  \:^{\eta_{\sigma}}_{p}\overline{X}^{(i)}_{n}(\omega),\ \forall(\omega.\sigma) \in \Omega
  \times \mathcal{S}^{(i)}$, \newline
 where $ \eta_{\sigma} \in \{0,1\}^{\mathcal{V}_{0}} $ is such that
   \small $ \eta_{\sigma}(2m) \stackrel{\mathrm{def}}{=} \left\{
  \begin{array}{clr}
  0 & , & m>0 \\
  \sigma (-2m) & , & -1 \leq m \leq -i \\
  1 & , & m<-i
  \end{array} \right. $;
 \normalsize
so that the IMC 
 $\:^{\eta_{\sigma}}_{p}Y^{(i)}_{\bullet}$ starts on state $ \sigma \in \mathcal{S}^{(i)}$. 
\end{itemize}

Now we state
 \begin{lem}
\label{lema7}
$ \mathbb E \left[ \frac{\:_{p}\overline{X}^{(i)}_{n}}{n}\right] \geq
  M^{(i)}(p), \:\forall n,i \in \mathbb N,\:
  p \in (0,1] $.
\end{lem}
\small
\textit{Proof:} \par
\begin{figure}[h]
\hspace{4cm}
\scalebox{.4}{\includegraphics[2.0cm,2.0cm][22.5cm,17.5cm]{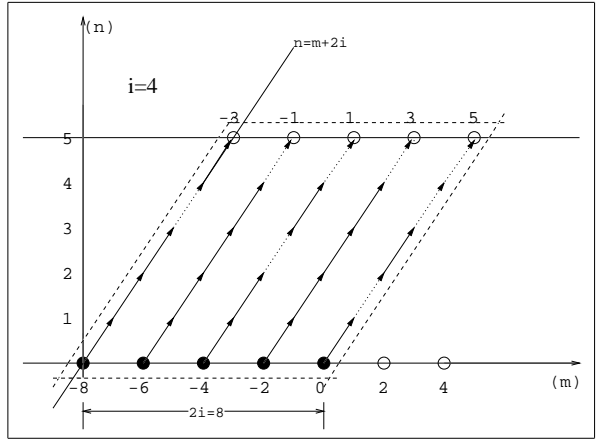}}
\vspace{1cm}
\caption{\label{figura3} \small No reinfection can take place in the region inside the
  dashed parallelogram.  The processes $\:_{p}X^{(4)}_{\bullet}$ and
  $\:_{p}X^{(\infty)}_{\bullet}$ are indistinguishable there.\normalsize     }
\end{figure}

 Inequality~\ref{acopl3} implies that
 \(
   \forall ( \omega , \sigma ) \in \Omega \times
   \mathcal{S}^{(i)}  , 
   \:_{p}\overline{\mathcal{X}}^{(i)}_{n} \left(
     (\omega , \sigma ) \right) = \:_{p}\overline{X}^{(i)}_{n}( \omega ) \geq
 \:^{ \eta_{ \sigma } }_{p} \overline{X}^{(i)}_{n} ( \omega ) =
 \:^{\pi^{(i)}}_{p} \overline{\mathcal{X}}^{(i)}_{n} \left(
    ( \omega , \sigma ) \right)  
\).
Integrating \mbox{ ( with respect to $ \mathbb P \times \pi^{(i)}$)}  both sides of
the inequality above, yields:
\[ \mathbb E \left[ \:_{p} \overline{\mathcal{X}}^{(i)}_{n} \right] = \left( 
    \mathbb E_{ \mathbb P } \left[ \:_{p} \overline{X}^{(i)}_{n} \right] \right) \geq 
 \mathbb E \left[  \:^{\pi^{(i)}}_{p} \overline{\mathcal{X}}^{(i)}_{n} \right]
\left( = nM^{(i)}(p) \right),  \]
where  the last  equality above  comes from the fact that the process   $\left\{ \:^{\pi^{(i)}}_{p}\overline{\mathcal{X}}^{(i)}_{n} \right\}_{n \in \mathbb
    N}$ has stationary increments of mean $M^{(i)}(p)$.
\begin{flushright} $\Box$ \end{flushright}
\normalsize
Figure~\ref{figura3} below refers to the next lemma :\vspace{-1cm} \newline
 
\begin{lem}
\label{lema8}
$ \lim_{i \rightarrow \infty} \mathbb E \left[ \:_{p}\overline{X}^{(i)}_{n}
\right] = \mathbb E \left[ \:_{p}\overline{X}^{(\infty)}_{n}
\right]$ for all (fixed) $ n \in \mathbb N $ and $ p \in (0,1]$.
\end{lem}
\small
\textit{Proof :} \par
 As the jumps to the right are bounded by $+1$, it follows that $\forall
 \omega \in \Omega , $
 \mbox{\( \:_{p}\overline{X}^{(i)}_{n} \leq n  
 \Longrightarrow  \:_{p}\overline{X}^{(i)}_{n} - m \leq n - m   . \) }
However, for a site $ (n,m) $ to be infected by force, we must have $ 
 \:_{p}\overline{X}^{(i)}_{n} - m > 2i $ (Def.~\ref{dfrec}(b) above). Hence,\linebreak
 if  $ n -m \leq 2i $, then $ \:_{p}\overline{X}^{(i)}_{n} - m \leq  2i,
  \forall \omega \in \Omega $ and site (n,m) will  not be infected \emph{ by
    force}, for any $ \omega \in \Omega $ 
. Accordingly, it can be proved by induction on $n$ that `` $  n-m \leq 2i \Longrightarrow X^{(\infty)}_{n}(m)[\omega] =  X^{(i)}_{n}(m)[\omega],
 \forall \omega \in \Omega $'', ie  all  sites $(n,m)$ such that $ n-m \leq 2i $
 have the same infection state regarding the processes $ \:_{p}X^{(\infty)}_{\bullet} $ and $
\:_{p}X^{(i)}_{\bullet} $. ( See Figure~\ref{figura3} above for the case: $i=4$.) 
\par Hence, for any (fixed) time $n$,
\begin{equation}
\label{Belitsky}
\:_{p}\overline{X}^{(\infty)}_{n} < \:_{p}\overline{X}^{(i)}_{n}
\Longrightarrow \:_{p}X^{(\infty)}_{n}(m)=0, \forall m \geq n-2i 
\end{equation}
The reasoning behind (~\ref{Belitsky}) is as follows:
\[ \:_{p}\overline{X}^{(i)}_{n} \geq n -2i \Longrightarrow
\:_{p}X^{(\infty)}_{n}\left(\:_{p}\overline{X}^{(i)}_{n}\right)=
\:_{p}X^{(i)}_{n}\left(\:_{p}\overline{X}^{(i)}_{n}\right)\stackrel{\mathrm{def}}{=}1\stackrel{(~\ref{acopl1})}{\Longrightarrow}\]
\[ \:_{p}\overline{X}^{(i)}_{n}\stackrel{\mathrm{def}}{=}\sup\{m \in \mathbb
Z_{n}:\:_{p}X^{(i)}_{n}(m)=1\}=\sup\{m \in \mathbb
Z_{n}:\:_{p}X^{(\infty)}_{n}(m)=1\}\stackrel{\mathrm{def}}{=}
\:_{p}\overline{X}^{(\infty)}_{n}\Longrightarrow  \]
\[\left\{ \:_{p}\overline{X}^{(i)}_{n}\geq n-2i\right\}\subset\left\{ \:_{p}\overline{X}^{(i)}_{n}= \:_{p}\overline{X}^{(\infty)}_{n}\right\}\Longrightarrow
\left\{ \:_{p}\overline{X}^{(i)}_{n}\neq
  \:_{p}\overline{X}^{(\infty)}_{n}\right\}\subset\left\{
  \:_{p}\overline{X}^{(i)}_{n} <  n-2i\right\}\subset\left\{
  \:_{p}\overline{X}^{(\infty)}_{n} <  n-2i\right\}\]
(~\ref{Belitsky}) follows identifying $\left\{ \:_{p}\overline{X}^{(i)}_{n}\neq
  \:_{p}\overline{X}^{(\infty)}_{n}\right\}$ with $\left\{
  \:_{p}\overline{X}^{(i)}_{n} >
  \:_{p}\overline{X}^{(\infty)}_{n}\right\}$ and $\left\{
  \:_{p}\overline{X}^{(\infty)}_{n} <  n-2i\right\}$ with \linebreak 
 $ \left\{
  \:_{p}X^{(\infty)}_{n}(m)=0, \forall m \geq n-2i \right\}$. 
\par Now, observe that, if the infection is not present on the set $\{(n,m)\in
  \mathcal{V}_{n} : m \geq n-2i \}$, all the paths joining its sites to slice
  $ \mathcal{V}_{0}$ must be obstructed somewhere. In particular all  the
  (i+1) straight lines joining site $ (0,-2j)$ to site $(n,n-2j),\: 0
\leq j \leq i \ $  must be interupted at some point (Figure ~\ref{figura3}
  above). As these lines 
  are made of different, independent bonds, the probability of this event equals
  $(1-p^{n})^{i+1}$, and we have
\[
\mathbb P \left( \:_{p} \overline{X}^{(i)}_{n} \neq
  \:_{p} \overline{X}^{(\infty)}_{n} \right) 
\leq \mathbb P \left( \:_{p}X^{(\infty)}_{n}(m)=0, \forall m \geq n-2i \right)
\leq (1-p^{n})^{i+1}  \]
So that $\:_{p}\overline{X}^{(i)}_{n}
  \stackrel{i}{\longrightarrow}  \:_{p}\overline{X}^{(\infty)}_{n}$ in
  probability, and there is a sub-sequence $(i_{k})_{k \in \mathbb N}$ such
  that $\:_{p}\overline{X}^{(i_{k})}_{n}
  \stackrel{k}{\longrightarrow}  \:_{p}\overline{X}^{(\infty)}_{n} \mathit {a\
  s}$ \linebreak~\cite[Theorem 7.6]{Bartle}. As the whole sequence
  $\left\{\:_{p}\overline{X}^{(i)}_{n}\right\}_{i \in \mathbb N}$ is non
  increasing (in i), we must have $\:_{p}\overline{X}^{(i)}_{n}
  \stackrel{i}{\longrightarrow}  \:_{p}\overline{X}^{(\infty)}_{n}\
  \mathit{a\:s\:}$ as well, and we can apply the \textit{Monotone Convergence Theorem} to
  conclude that $\mathbb E \left[ \:_{p}\overline{X}^{(i)}_{n}\right]
  \searrow^{i}\mathbb E \left[
  \:_{p}\overline{X}^{(\infty)}_{n}\right],\:\forall n \in \mathbb N$.

\begin{flushright}$ \Box $ \end{flushright}
\normalsize
Now we can prove the \emph{Convergence Theorem}:
\begin{teorema}
\label{conv.teo}
\( p^{(i)}_{c} \nearrow_{i} p_{c} \)
\end{teorema}
\small
\textit{Proof :}
Suppose that $\lim_{i \rightarrow \infty } p^{(i)}_{c} < p_{c}$, so that we can
choose $p$ such that $ \lim_{i \rightarrow \infty } p^{(i)}_{c}<p <p_{c}$. Then Lemma~\ref{lema2} and Corollary~\ref{coro5} ensure that $ \forall i \in \mathbb N
,\ M^{(i)}(p) > M^{(i)}(p^{(i)}_{c}) \stackrel{\mathrm{def}}{=}0 $ and Lemma~\ref{lema7} that $ \mathbb E \left[ \:_{p}\overline{X}^{(i)}_{n} \right] \geq
nM^{(i)}(p)>0 , \forall i \in \mathbb N,\forall n \in \mathbb N $. Thus
\begin{equation}
\label{absurd}
 \lim_{i \rightarrow \infty } \mathbb E \left[ \:_{p}\overline{X}^{(i)}_{n}
\right]\geq 0,\forall n \in \mathbb N
 \end{equation}
(inequality~\ref{acopl1} ensures that this limit is well defined)
\par By the other side,
\begin{eqnarray*}
p<p_{c} & \Longrightarrow \frac{\:_{p}\overline{X}^{(\infty)}_{n}}{n}
\stackrel{(n)}{\longrightarrow} - \infty \ \mathit{a\ s\ } & \mbox{
  Lemma~\ref{lema1} and commentary on Lemma~\ref{lema2}} 
\\ & \Longrightarrow \mathbb E \left[
  \frac{\:_p\overline{X}^{\infty}_{n}}{n} \right] \stackrel{(n)}{\longrightarrow} -
\infty & \mbox{ Fatou's Lemma} \\
 & \Longrightarrow \exists \overline{n} : \mathbb E \left[
  \:_p\overline{X}^{\infty}_{\overline{n}} \right] <0 & \end{eqnarray*}
Applying Lemma~\ref{lema8} for $\overline{n}$ yields
\[ \lim_{i\rightarrow\infty} \mathbb E \left[
  \:_p\overline{X}^{i}_{\overline{n}} \right]= \mathbb E \left[
  \:_p\overline{X}^{\infty}_{\overline{n}} \right]<0 \]

which contradicts inequality~\ref{absurd} above. Hence,  $\lim_{i
  \rightarrow \infty } p^{(i)}_{c} \geq  p_{c}$ and the theorem follows from
  Corollary~\ref{coro6}.
\begin{flushright} $\Box$ \end{flushright}
\normalsize

\section{Numerical Calculations}
\label{CalcNum}
\subsection{Algebraical Determination of the Critical Probabilities}
\label{algoritmo}
Theorem~\ref{conv.teo} (\emph{The Convergence Theorem}) is of theoretical interest by itself. However weaker results such as Corollaries~\ref{coro5} and~\ref{coro6}  already indicate that each $p^{(i)}_{c},\: i \in \mathbb N $, is an improved lower bound (regarding its predecessors) for the critical probability of \emph{Oriented Percolation}.\par
Although tacitly present in Sections~\ref{CMIs} and~\ref{probcrit} above, we present below the  algorithm for calculating the \emph{critical probabilities} $p^{(i)}_{c},\: i \in \mathbb N $ exactly.\par
According to Definitions ~\ref{IMC},~\ref{Pilmk},~\ref{Mil},~\ref{Mdrift} and ~\ref{dfpc}, the \emph{critical probabilities} $p^{(i)}_{c},\: i \in \mathbb N $, may be determined in the following steps:
\begin{description}
\item[(i)] \textbf{Determination of the Probabilities $\mathit{ P^{(i)}_{(l,m,k)}} $} (Def.~\ref{Pilmk}) \textbf{in terms of} $\mathit{ q \stackrel{\mathrm{def}}{=}1-p:} $ \par
elementary combinatorics show that the probabilities $ P^{(i)}_{(l,m,k)},0\leq l,m \leq 2^{i}-1,k \in \mathbb N $, may be expressed as polinomial functions of $q$, ie
\begin{equation} \label{Pilmkq} P^{(i)}_{(l,m,k)}= P^{(i)}_{(l,m,k)} (q) : \mbox{polinomial in } \:q \end{equation}
\item[(ii)] \textbf{Determination of the Transition Matrix  $\left( q^{(i)}_{lm}\right)_{0\leq l,m\leq 2^{i}-1}$ of the IMC} $\mathit{\:_{p}Y^{(i)}_{\bullet}: }$ \par
the transition probabilities defined in (~\ref{trans.prob}) may be expressed as 
\begin{equation}
\label{qSP}
q^{(i)}_{lm}=\sum^{\infty}_{k=0} P^{(i)}_{(l,m,k)} (q)
\end{equation}
As the numerical sequence $ \left( P^{(i)}_{(l,m,k)} (q)  \right)_{k\in\mathbb N}$ is a \emph{geometric progression} (ratio $q^2$) starting from the $(i+2)^{th}$ term, the transition probabillities $q^{(i)}_{lm}$ may be expressed as rational functions of $q$. As a matter of fact, these probabilities are always polinomials in $q$:
\[q^{(i)}_{lm}=q^{(i)}_{lm}(q) : \mbox{polinomial in} \:q \]
\item[(iii)] \textbf{Determination of the Stationary Measure } $\mathit{(\pi^{i}) }$ \textbf{of the IMC} $\mathit{\:_{p}Y^{(i)}_{\bullet}: }$\par
the combinatorial calculus that lead to equation (~\ref{Pilmkq}) above show that the \emph{transition probabilities} $ q^{(i)}_{lm}(q) ,0\leq l,m\leq 2^{i}-1 $, are strictly positive for $q\in (0,1)$, so that the \emph{transition matrix} $\left( q^{(i)}_{lm}\right)$ is irreducible and aperiodic. Hence $\pi^{i}$ is a well defined probability measure on $\mathcal{S}^{i}$ and can be algebraically determined from  $\left( q^{(i)}_{lm} \right)$; resulting that
\[\pi^{(i)}_{l}(q) : \mbox{\emph{a strictly positive rational function of}} \:q,\:q\in (0,1),\:0\leq l\leq 2^{i}-1 \]
\item[(iv)] \textbf{Determination of } $\mathit{M^{(i)}_{l}}$, \textbf{the Mean Jump of the SDTCP} $\mathit{\:_{p}X^{(i)}_{\bullet}}$ \textbf{on State} $\mathit{\sigma_{l}}:$
\par
Definition~\ref{Mil} and the fact that the numerical sequence $\left( P^{(i)}_{lmk}(q)\right)_{k\in\mathbb N}$ is essentially a \emph{geometric progression} yield that
\[M^{(i)}_{l}(q) : \mbox{ \emph{a rational function of} } \:q,\ 0\leq l\leq 2^{i}-1\]
\item[(v)] \textbf{Determination of } $\mathit{M^{(i)}}$, \textbf{ the Mean Jump of the SDTCP } $\mathit{\:_{p}X^{(i)}_{\bullet}:}$ \par
Definition~\ref{Mdrift} and steps (iii) and (iv) above ensure that 
\[M^{(i)}(q) : \mbox{ \emph{ a rational function of }} q \]
\item[(vi)] \textbf{ Determination of the Critical Probability } $ \mathit{p^{(i)}_{c} } :$ \par
according to Corollary~\ref{coro3} and Definition~\ref{dfpc},  $ \mathit{p^{(i)}_{c} } $ is the only real root of $M^{(i)}(p)$ on the interval $(0,1)$; so that it is the only real root of the polinomial in the numerator of $M^{(i)}(p)$.  Its calculation thus can be done numerically.    
\end{description}
In the sequel, we employ the algorithm described above for calculating the first critical probability in algebraical terms:
\subsubsection{The Critical Probability of Zero $\!\!\!\!\,^{th}$ Order , $ p^{(0)}_{c}$ :} 
In this case $ \left| \mathcal{S}^{(i)} \right| =1$, so $ \:_{p}Y^{(i)}_{\bullet} $ is trivial:
\begin{eqnarray*}
M^{(0)}(q) & = & M^{(0)}_{0}(q) = \sum^{\infty}_{k=0} (1-2k).P^{(0)}_{0,0,k}(q)=1.p-1.q(1-q^{2})-3.q^{3}(1-q^{2})-5.q^{5}(1-q^{2}) - \ldots \\
           & = & (1-q) -q(1-q^{2}).\left[1+3q^{2}+5q^{4}+7q^{6}+\ldots \right]=(1-q)-q(1-q^{2})\frac{1+q^{2}}{(1-q^{2})^{2}}=\\
           & = & (1-q) - \frac{q+q^{3}}{1-q^{2}}=\frac{1-2q-q^{2}}{1-q^{2}}
\end{eqnarray*}
Hence, 
\[ M^{(0)}(q)=0\Leftrightarrow 1-2q-q^{2}=0\Rightarrow   p^{(0)}_{c}=2-\sqrt{2}=0.58579\ldots \]
\subsubsection{The Critical Probability of First Order, $p^{(1)}_{c}$ :}
\[ \left( q^{(1)}_{lm}\right)(q)=   \left[   \begin{array}{cc}
                                         q-q^{3}+q^{4} & 1-q+q^{3}-q^{4}\\
                                         q^{2}         & 1-q^{2}           \end{array} \right], \ \ \ \ \pi^{(1)}_{l}(q)= \frac{\left(q^{2},1-q+q^{3}-q^{4} \right)}{1-q+q^{2}+q^{3}-q^{4}} \]
\[M^{(1)}_{0}(q)=\frac{1-2q-3q^{2}+2q^{4}}{1-q^{2}},\ \ M^{(1)}_{1}(q)=\frac{1-2q-q^{2}}{1-q^{2}}\]
\[M^{(1)}(q)=\frac{1-3q+2q^{2}-6q^{4}+q^{5}+3q^{6}}{1-q+2q^{3}-2q^{4}-q^{5}+q^{6}} \]
Hence,
\[M^{(1)}(q)=0 \Rightarrow p^{(1)}_{c}=0.604233\ldots\]
\subsubsection{The Critical Probability of Second Order, $p^{(2)}_{c}$ :}
\scriptsize

\[ \left( q^{(2)}_{lm}(q) \right) =      \left[ \begin{array}{cccc}
          q-q^{3}+q^{6} & q-2q^{2}+q^{3}+q^{4}-q^{6} & 1-2q+q^{2}+q^{4}-q^{6} & q^{2}-2q^{4}+q^{6} \\
          q^{2}-q^{3}+q^{4}-q^{5}+q^{6} & q-q^{2}+q^{5}-q^{6} & q-2q^{2}+2q^{3}-q^{4}+q^{5}-q^{6} & 1-2q+2q^{2}-q^{3}-q^{5}+q^{6} \\
          2q^{3}-q^{4}-q^{5}+q^{6}-q^{7}+q^{8} & 2q^{2}-3q^{3}+q^{4}+q^{7}-q^{8} & q-2q^{3}+2q^{5}-q^{6}+q^{7}-q^{8} & 1-q-2q^{2}+3q^{3}-q^{5}-q^{7}-q^{8} \\
          q^{4} & q^{2}-q^{4} & q^{2}-q^{4} & 1-2q^{2}+q^{4} 
          \end{array} \right] \]
\[ \pi^{(2)}_{0}(q) = \frac{  -2q^{4}+2q^{5}-4q^{6}+9q^{7}-14q^{8}+15q^{9}-9q^{10}+2q^{11}}{-1+3q-6q^{2}+8q^{3}-17q^{4}+30q^{5}-44q^{6}+46q^{7}-20q^{8}-17q^{9}+38q^{10}-32q^{11}+13q^{12}-2q^{13}} \]
\[ \pi^{(2)}_{1}(q) = \frac{ -q^{2}+2q^{3}-2q^{4}+q^{5}-4q^{6}+9q^{7}-5q^{8}-6q^{9}+12q^{10}-8q^{11}+2q^{12}}{-1+3q-6q^{2}+8q^{3}-17q^{4}+30q^{5}-44q^{6}+46q^{7}-20q^{8}-17q^{9}+38q^{10}-32q^{11}+13q^{12}-2q^{13}} \]
\[ \pi^{(2)}_{2}(q) = \frac{-q^{2}+q^{3}+2q^{5}-9q^{6}+15q^{7}-14q^{8}+8q^{9}-2q^{10} }{-1+3q-6q^{2}+8q^{3}-17q^{4}+30q^{5}-44q^{6}+46q^{7}-20q^{8}-17q^{9}+38q^{10}-32q^{11}+13q^{12}-2q^{13}}  \]
\[ \pi^{(2)}_{3}(q) = \frac{-1+3q-4q^{2}+5q^{3}-13q^{4}+25q^{5}-27q^{6}+13q^{7}+13q^{8}-34q^{9}+37q^{10}-26q^{11}+11q^{12}-2q^{13}}{-1+3q-6q^{2}+8q^{3}-17q^{4}+30q^{5}-44q^{6}+46q^{7}-20q^{8}-17q^{9}+38q^{10}-32q^{11}+13q^{12}-2q^{13}} \]
\[M^{(2)}_{l} = \frac{1}{1-q^{2}} \left(1-2q-5q^{2}+4q^{4},1-2q-3q^{2}+2q^{4},1-2q-q^{2}-2q^{4}+2q^{6},1-2q-q^{2}\right)\]
\[ M^{(2)}(q) = \frac{-1+5q-11q^{2}+16q^{3}-25q^{4}+52q^{5}-75q^{6}+96q^{7}-58q^{8}-69q^{9}+152q^{10}-111q^{11}-5q^{12}+74q^{13}-49q^{14}+10q^{15} }{\left(1-q^{2}\right)\left(-1+3q-6q^{2}+8q^{3}-17q^{4}+30q^{5}-44q^{6}+46q^{7}-28q^{8}-17q^{9}+38q^{10}-32q^{11}+13^{12}-2q^{13} \right) } \]
\small
Hence,
\[M^{(2)}(q)=0 \Rightarrow p^{(2)}_{c}=0.614187\]
\normalsize
\subsection{Numerical Determination of the Critical Probabilities:}
The combinatorial calculations leading to the transition matrix $\left(q^{(i)}\right)_{lm}$ and to the mean-drift vector $M{(i)}_{l}$ described in Section ~\ref{algoritmo} above, yet fastidious for humans, are taylor-made for computers: while it took us a whole afternoom for calculating  $\left(q^{(2)}\right)_{lm}$ algebraically (16 polinomial entries), a FORTRAN 77 program, running at 600 MHz, calculated  $\left(q^{(9)}\right)_{lm}$ ($ 4^9$ polinomial entries) in less than one minute time! \par
The greatest problem in writing a computer program version for the algorithm described in Section ~\ref{algoritmo} arises precisely in step(iii), viz. the determination of the stationary measure $\pi^{(i)}$ (in algebraical terms), starting from the transition matrix  $\left(q^{(i)}\right)_{lm}(q)$. True, it can be done straightforwardly, acting on the polinomials as if they were numbers. However, if we don't take into account the eventual and probable simplifications that may take place, the degree of the outcoming polinomials soon get unmanageble. And finding out these simplifications is far from obvious.\par
Therefore we have adopted a different approach:
\begin{description}
\item[(i)] After having algebraically calculated the transition matrix  $\left(q^{(i)}\right)_{lm}(q)$ and the mean-drift vector $M^{(i)}_{l}(q)$, the program calculated a numerical transition matrix and a numerical mean-drift vector for a decreasing sequence of numerical values of $q$.
\item[(ii)] From each numerical transition matrix, a numerical stationary
  measure was obtained by solving a system of linear equations. It is
  worth-while mentioning that \emph{partial pivoting}~\cite{Wilkinson} was of
  fundamental importance in keeping numerical errors under controll.
\item[(iii)] The numerical \emph{mean-drift} was then obtained, performing the inner-product of Definition~\ref{Mdrift}.
\item[(iv)] According to Definition~\ref{dfpc} $\ q^{(i)}_{c} (\stackrel{\mathrm{def}}{=}1-p^{(i)}_{c})$ lies between the last value of $q$ for which $M^{(i)}(q)$ is negative and the first value of $q$ for which $M^{(i)}(q)$ is positive. 
\end{description}
This approach has the onus of introducing numerical erros (basically in step ii above), from which we hitherto have got rid of. So as to produce reliable numerical data for the critical probabilities, rigorous  upper bounds for the numerical errors must be provided. In doing so, we followed the set-up  of \emph{foward analysis} described in ~\cite{Wilkinson}, which consists fundamentally in providing partial upper bounds after each arithmetic operation performed in the flow of numerical calculations leading to the final numerical result. This rather crude approach has the advantage of providing a secure upper bound for  numerical errors, regardless of particular features the linear system may present, as it is often the case in the so called \emph{backward analysis} (~\cite{Wilkinson}).\par
In the sequel, we present the numerical data produced  following\vspace{11 mm}
steps i-iv above\footnote{the calculations were carried out on double precision. } :\par\par\par
\small
\hspace{-7 mm}
\mbox{ \begin{tabular}{|c|c|c|}
\multicolumn{3}{c}{i=1}\\ \hline
            p & $M^{(1)}_{Num}(p)$ & $\Delta M^{(1)}(p)$ \\ \hline

           0.604231 & -9.15E-6 & 6.2E-16 \\
           0.604232 & -5.27E-6 & 6.2E-16 \\
           0.604233 & -1.40E-6 & 6.2E-16 \\ 
           0.604234 & +2.47E-6 & 6.2E-16 \\ \hline
\end{tabular}\hspace{3 mm}  \begin{tabular}{|c|c|c|}
\multicolumn{3}{c}{i=2}\\ \hline
            p & $M^{(2)}_{Num}(p)$ & $\Delta M^{(2)}(p)$ \\ \hline

           0.614185 & -1.01E-5 & 6.1E-15 \\
           0.614186 & -5.85E-6 & 6.2E-15 \\
           0.614187 & -1.57E-6 & 6.1E-15 \\ 
           0.614188 & +2.71E-6 & 6.2E-15 \\ \hline
\end{tabular}\hspace{3 mm} \begin{tabular}{|c|c|c|}
\multicolumn{3}{c}{i=3}\\ \hline
            p & $M^{(3)}_{Num}(p)$ & $\Delta M^{(3)}(p)$ \\ \hline
         
           0.620203 & -1.36E-5 & 2.5E-14 \\
           0.620204 & -8.96E-6 & 2.5E-14 \\
           0.620205 & -4.32E-6 & 2.6E-14 \\
           0.620206 & +3.15E-7 & 2.6E-14 \\ \hline

\end{tabular} }\vspace{10 mm}\par  
    
\hspace{-7 mm} \vspace{7 mm}
\mbox{ \begin{tabular}{|c|c|c|}
\multicolumn{3}{c}{i=4}\\ \hline
            p & $M^{(4)}_{Num}(p)$ & $\Delta M^{(4)}(p)$ \\ \hline

           0.624210 & -9.52E-6 & 8.0E-14 \\
           0.624211 & -4.57E-6 & 8.0E-14 \\
           0.624212 & +3.71E-7 & 8.0E-14 \\ 
                    &          &         \\ \hline
\end{tabular} \hspace{3 mm} 
                                                                                                                                                                                                                \begin{tabular}{|c|c|c|}
\multicolumn{3}{c}{i=5}\\ \hline
            p & $M^{(5)}_{Num}(p)$ & $\Delta M^{(5)}(p)$ \\ \hline
         
           0.627064 & -1.31E-5 & 1.9E-13 \\
           0.627065 & -7.85E-6 & 1.9E-13 \\
           0.627066 & -2.63E-6 & 1.9E-13 \\
           0.627067 & +2.60E-6 & 1.9E-13 \\ \hline

\end{tabular}  \hspace{3 mm}                                                       

 \begin{tabular}{|c|c|c|}
\multicolumn{3}{c}{i=6}\\ \hline
            p & $M^{(6)}_{Num}(p)$ & $\Delta M^{(6)}(p)$ \\ \hline
         
           0.629201 & -1.55E-5 & 3.8E-13 \\
           0.629202 & -9.98E-6 & 3.8E-13 \\
           0.629203 & -4.51E-6 & 3.9E-13 \\
           0.629204 & +9.66E-7 & 3.8E-13 \\ \hline

\end{tabular} }
\hspace{-7 mm} \vspace{7 mm}
\mbox{ \begin{tabular}{|c|c|c|}
\multicolumn{3}{c}{i=7}\\ \hline
            p & $M^{(7)}_{Num}(p)$ & $\Delta M^{(7)}(p)$ \\ \hline

           0.630863 & -1.01E-5 & 7.2E-13 \\
           0.630864 & -4.41E-6 & 7.2E-13 \\
           0.630865 & +1.30E-6 & 7.2E-13 \\ 
                    &          &         \\ \hline
\end{tabular} \hspace{3 mm}
 \begin{tabular}{|c|c|c|}
\multicolumn{3}{c}{i=8}\\ \hline
            p & $M^{(8)}_{Num}(p)$ & $\Delta M^{(8)}(p)$ \\ \hline

           0.632192 & -8.53E-6 & 1.3E-12 \\
           0.632193 & -2.60E-6 & 1.2E-12 \\
           0.632194 & +3.32E-6 & 1.3E-12 \\ 
                    &          &         \\ \hline
\end{tabular} 
\hspace{3 mm}
  \begin{tabular}{|c|c|c|}
\multicolumn{3}{c}{i=9}\\ \hline
            p & $M^{(9)}_{Num}(p)$ & $\Delta M^{(9)}(p)$ \\ \hline
         
           0.63326 & -1.29E-4 & 2.0E-12 \\
           0.63327 & -6.79E-5 & 2.0E-12 \\
           0.63328 & -6.68E-6 & 2.0E-12 \\
           0.63329 & +5.45E-5 & 2.0E-12 \\ \hline
 
\end{tabular} }\footnote{$|M^{(i)}(p)-M^{(i)}_{Num}(p)|\leq\ \Delta M^{(i)}(p)$
           ; where $ M^{(i)}_{Num}(p)$ denotes the numerically calculated
             \emph{mean drift} of $ i^{th}$ order.}
\normalsize
 
Based upon the data displayed above, the first ten critical probabilities may
             be determined:\vspace{4 mm}\small 
 
\hspace{-7 mm} \mbox{ \begin{tabular}{|c||c|c|c|c|c|c|c|c|c|c|} \hline
                $ i$ & 0 & 1 & 2 & 3 & 4 & 5 & 6 & 7 & 8 & 9 \\ \hline 
 $p^{(i)}_{c}$ & 0.585787 & 0.604233 & 0.614187 & 0.620205 & 0.624211 &
0.627066 & 0.629203 & 0.630864 & 0.632193 & 0.63328 \\ \hline
\end{tabular}}\vspace{10 mm}\normalsize\par
So  that we can establish the inequality : \hspace{10 mm} \vspace{17 mm}\large \fbox{ $ p_{c} \geq 0.63328 $} 
\normalsize
\section{Simulations}
The boolean and inductive features of definition ~\ref{dfrec} make it
extremely suitable for computer simulations: at each step a FORTRAN77-written
computer program determmines the value of $ \:_{p}\overline{X}^{(i)}_{n}/n $, the right
edge mean speed at moment $n$. In accordance with lemma ~\ref{lema1} this
sequence converges ($\mathit{a.s.}$  in $n$) to $M^{(i)}(p)$, the right edge mean
speed. Therefore running the program for different values of $p$, one can
estimate the value of $p_{c}^{(i)}$ observing the height of the plateau
established with the increase of $n$, according with the rule:
\begin{description}
\item[(a)] \emph{height of plateau} $<0 \Rightarrow p<p^{(i)}_{c}$,
\item[(b)] \emph{height of plateau} $>0 \Rightarrow p>p^{(i)}_{c}$.
\end{description}
The two figures displayed below ilustrate the use of this technique when
$i=1000$:
\begin{figure}[h]
\hspace{2cm}
\scalebox{.40}[.25]{ \rotatebox{-90}{\includegraphics[20cm,10cm]{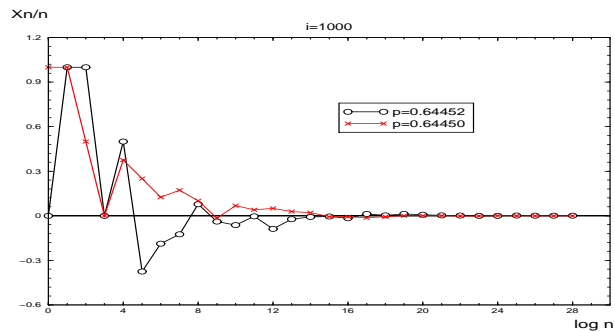}}}
\caption{\label{Grafi1000a}Independent (not coupled) trajectories of the processes
  $\:_{.64450}\overline{X}^{(1000)}_{n}/n$ e
  $\:_{.64452}\overline{X}^{(1000)}_{n}/n$.}
\end{figure}

\begin{figure}[h]
\hspace{2cm}
\scalebox{.40}[.25]{ \rotatebox{-90}{\includegraphics[20cm,10cm]{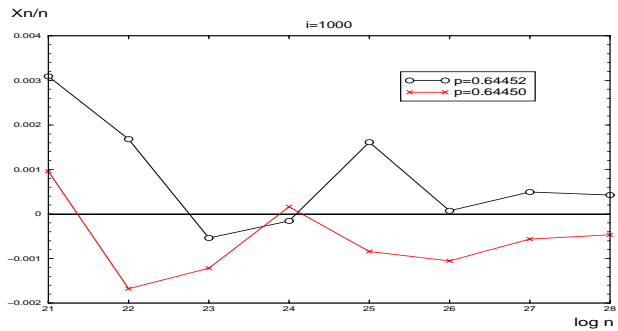}}}
\caption{Enlarged detail of the figure above .}
\end{figure}
\newpage 
Table ~\ref{tabsimul} below sumarizes the numerical results obtained by means
of the technique for different values of $i$ :
\begin{table}[h]
\large
\begin{center}
\begin{tabular}{||c|c||} \hline \hline
                  $i$ & $p^{(i)}_{c}$ \\ \hline
                   5 &    0,627 \\
                   6 &    0,629 \\
                   9 &    0,6332 \\
                   20 &   0,638 \\
                   40 &   0,641 \\
                   100 &  0,643 \\
                   200 &  0,6438 \\
                   1000 & 0,64451 \\ \hline \hline 

\end{tabular}
\end{center}
\normalsize
\caption{\label{tabsimul}Critical probabilities obtained by means of simulation. }
\end{table}
\\
It is interesting to observe that the \emph{plateau pattern} exhibited in
figure ~\ref{Grafi1000a} is present even in the vicinity of criticality,
\emph{ie} $ p \approx p^{(i)}_{c}$ and $ i \nearrow \infty $.
Therefore the critical probabilities $p^{(i)}_{c}$ (even for large values of
$i$ and so very close to criticality) can be estimated within the desired
precision simply by increasing the value of $n$, \emph{ie} running the program
until the \emph{plateau pattern} becomes clear.
We belive that this feature distinguishes our technique from the usual
\emph{Monte Carlo} simulation methods, wherein instability appears near
criticality  and the use of scaling techniques is called for.
\\
\par \hspace{-6 mm} \Large \textbf{Acknowledgement:}\vspace{5 mm}\par
\normalsize The authors are indebted to M. Menshikov, who first suggested the
use of\emph{ Markov chains in a strip} to \emph{track} the edge of Oriented Percolation\ldots


\begin{thebibliography}{99}
\bibitem{Durrett_1} R. Durrett. (1984)
        \emph{Oriented Percolation in Two Dimensions.}\par
        The Annals of Probability, 12(4), 999-1040.
\bibitem{Grimmett_1} G. Grimmett. (1999)
        \emph{Percolation.} Springer Verlag.
\bibitem{Misha}G.Fayolle/V.A.Malyshev/M.V.Menshikov. (1995) \par
        \emph{Topics in the Constructive Theory of Countable Markov Chains}. 
        Cambridge University Press .
\bibitem{Bartle} R. G. Bartle. (1995)\par \emph{The Elements of Integration
        and Lebesgue Measure.} John Willey \& Sons Inc.
\bibitem{Wilkinson} J. H. Wilkinson. (1963)
        \emph{Rounding Errors in Algebraic Processes.} Her Majesty's Stationary Office.
\end{thebibliography}
\end{document}